\documentclass[german,english
              ,twoside
              ,11pt
              ]{amsart}
\usepackage[T1]{fontenc}
\usepackage[pagesize,BCOR1cm,DIV11]{typearea}
\usepackage[latin1]{inputenc}
\usepackage{babel}
\usepackage{hyperref}
\usepackage{cas-math}
\usepackage{cas-paper}
\usepackage[arrow,cmtip,matrix]{xy}

\DeclareMathOperator{\conn}{conn}
\DeclareMathOperator{\ind}{ind_{\Z_2}}
\DeclareMathOperator{\coind}{coind_{\Z_2}}
\DeclareMathOperator{\hind}{cohom-ind_{\Z_2}}
\DeclareMathOperator{\Mon}{Mon}
\def\real{\abs}
\newcommand{\FHom}{\mathrm F\Hom}

\newenvironment{acknowledgements}{\subsection{Acknowledgements}}{}
\newcommand{\subclass}[1]{}
\begin{document}

\title{Graph colourings, spaces of edges and spaces of circuits}

\author{Carsten Schultz}
\address{Institut für Mathematik, MA 6-2\\
Technische Universität Berlin\\
D-10623 Berlin, \hbox{Germany}}
\email{carsten@codimi.de}
\date{June 2006}
\thanks{This research was partially supported by the
\foreignlanguage{german}{Deutsche Forschungsgemeinschaft} within the
European graduate program ``Combinatorics, Geometry, and Computation''
(GRK~588/2)
and partially by the Leibniz grant of G.~M.~Ziegler.
Some of the results were obtained during a stay at the 
Institut Mittag-Leffler.
}

\begin{abstract}
By Lovász' proof of the Kneser conjecture, the chromatic number of a
graph~$G$ is bounded from below by the index of the $\Z_2$-space
$\Hom(K_2,G)$ plus two.  We show that the cohomological index of
$\Hom(K_2,G)$ is also greater than the cohomological index of the
$\Z_2$-space $\Hom(C_{2r+1}, G)$ for $r\ge1$.  This gives a new and
simple proof of the strong form of the graph colouring theorem by
Babson and Kozlov, which had been conjectured by Lovász, and at the
same time shows that it never gives a stronger bound than can be
obtained by $\Hom(K_2, G)$.  The proof extends ideas introduced by 
\v Zivaljevi\'c in a previous elegant proof of a special case.
We then generalise the arguments and obtain conditions under
which corresponding results hold for other graphs in place of
$C_{2r+1}$.  This enables us to find an infinite family of test graphs
of chromatic number~$4$ among the Kneser graphs.

Our main new result is a description of the $\Z_2$-homotopy type of
the direct limit of the system of all the spaces $\Hom(C_{2r+1}, G)$
in terms of the $\Z_2$-homotopy type of $\Hom(K_2, G)$.  A corollary
is that the coindex of $\Hom(K_2, G)$ does not exceed the coindex of
$\Hom(C_{2r+1}, G)$ by more then one if $r$~is chosen sufficiently
large.  Thus the graph colouring bound in the theorem by
Babson~\&~Kozlov is also never weaker than that from Lovász' proof of
the Kneser conjecture.
\subclass{57M15; 05C15}
\end{abstract}

\maketitle

\section{Introduction}
%%%%%%%%%%%%%%%%%%%%%%
\subsection{Background}
As a means of proving Kneser's Conjecture, Lovász has shown that a
graph is not $k$-colourable if its neighbourhood complex is
$(k-2)$-connected.  Since the neighbourhood complex of $G$ is homotopy
equivalent to the cell complex $\Hom(K_2, G)$, which was introduced
later, this result can be stated as follows.  All necessary
definitions will be given in the next section.
\begin{thm}[Lovász~\cite{lovasz}]\label{thm:L}
Let $G$ be a graph.  Then \[\conn\Hom(K_2, G)\le\chi(G)-3.\]
\end{thm}
The reformulation in terms of the $\Hom$-complex made it natural to
ask if similar theorems would hold for graphs other than $K_2$.  In
particular, one might have hoped that
$\conn\Hom(T,G)\le\chi(G)-\chi(T)-1$ for all graphs $T$ and $G$,
provided that $\Hom(T, G)\ne\emptyset$.  A graph~$T$ such that this
holds for all graphs~$G$ is called a \emph{test graph}~\cite{babson-kozlov-i}.  
Hoory and
Linial have shown that not every graph is a test graph by giving an
example of a graph $T$ with $\conn\Hom(T,
K_{\chi(T)})\ge0$~\cite{linial-test-graph}.  Babson and Kozlov
succeeded in proving the following positive result that had been
conjectured by Lovász.
\begin{thm}[Babson \& Kozlov~\cite{babson-kozlov-ii}]\label{thm:BK}
Let $G$ be a graph.  Then \[\conn\Hom(C_{2r+1}, G)\le\chi(G)-4.\]
\end{thm}

The spaces $\Hom(K_2, G)$ and $\Hom(C_{2r+1}, G)$ are equipped with
free $\Z_2$-actions.  For such spaces there are several index functions
that assign to the space an integer measuring the complexity of the
action.  We give definitions in \ref{def:ind} and recall that
\begin{equation*}
\conn X+1\le\coind X\le\hind X\le\ind X
\end{equation*}
for all free $\Z_2$-spaces~$X$.  In both of the above theorems, $\conn
X$ is only used in the statement for convenience, since it does not
depend on the action.  In both cases, $\conn X$ can immediately be
replaced by $\coind X-1$.

In the case of the original Lovász criterion, the now usual proof
yields the following result.
\begin{thm}\label{thm:Lext}
Let $G$ be a graph.  Then $\ind\Hom(K_2, G)\le\chi(G)-2$.
\end{thm}

\begin{proof}
A colouring $c\colon G\to K_n$ induces a $\Z_2$-map
$\Hom(K_2,G)\to\Hom(K_2,K_n)$, and $\Hom(K_2,K_n)$ is
$\Z_2$-homeomorphic to the $(n-2)$-sphere with the antipodal map.  For
this homeomorphism see e.g.~\cite[4.2]{babson-kozlov-i},
\cite[Rem.~2.5]{c5}, or \prettyref{ex:CoddKn}.
\end{proof}
\begin{sloppypar}
We remark that the inequality $\coind X\le\ind X$ that is used to
obtain \prettyref{thm:L} from \prettyref{thm:Lext} is essentially the
Borsuk-Ulam Theorem.
\end{sloppypar}

For $C_{2r+1}$, Babson and Kozlov had proposed and partially proven
the following slightly stronger version of \prettyref{thm:BK}.
\begin{thm}\label{thm:S}
Let $G$ be a graph.  Then $\hind\Hom(C_{2r+1}, G)\le\chi(G)-3$.
\end{thm}
\subsection{Results}
The main result of the current work is the following.
\begin{thm}\label{thm:main-intro}
Let $G$ be a graph.  Then
\begin{equation*}
\colim_{r}\Hom(C_{2r+1}, G)
    \homot_{\Z_2}\Map_{\Z_2}(\Sphere^1_b,\Hom(K_2,G)).
\end{equation*}
\end{thm}
This result will be part of \prettyref{thm:homot-colim}.  The left
hand side of this homotopy equivalence is the colimit of a diagram
that will be defined in \prettyref{sec:colimit}.  The right hand side
is the space of all equivariant maps from $\Sphere^1$ equipped with
the antipodal map to $\Hom(K_2, G)$.  This space is made into a
$\Z_2$-space via a second $\Z_2$-action on $\Sphere^1$ that is a
reflection by a line through the origin in~$\R^2$; that $\Sphere^1$ is
equipped with these two actions is what the notation $\Sphere^1_b$
indicates (see~\prettyref{nota:Skb}).

The space $\Hom(K_2, G)$ can be thought of as the space of oriented
edges of~$G$, the $\Z_2$-action being orientation reversal.  The space
$\colim_{r}\Hom(C_{2r+1}, G)$ can be thought of as the space of
parametrized circuits, or closed paths, in~$G$ of arbitrary but odd
length, the $\Z_2$-action being the reversal of the direction of the
closed paths.  The theorem describes that the $\Z_2$-homotopy type of
the former determines the $\Z_2$-homotopy type of the latter.

The consequences of this homotopy equivalence for the graph colouring
theorems above is described by the following corollaries.  Since in
the body of this article they will be proved before the theorem, we
here give their derivations from it.  Readers who are not very
familiar with the methods are encouraged to skip them.
\begin{cor}\label{cor:ineq1}
Let $G$ be a graph with at least one edge.  Then
\[\hind\Hom(C_{2r+1},G)+1\le\hind\Hom(K_2,G).\]
\end{cor}
\begin{proof}
We have a composition of $\Z_2$-maps
\begin{align*}
\Sphere^1_b\times_{\Z_2}\Hom(C_{2r+1}, G)
&\to_{\Z_2}\Sphere^1_b\times_{\Z_2}\colim_r\Hom(C_{2r+1}, G)\\
&\homot_{\Z_2}\Sphere^1_b\times_{\Z_2}\Map_{\Z_2}(\Sphere^1_b,\Hom(K_2,G))\\
&\to_{\Z_2}\Hom(K_2,G),
\end{align*}
where the last arrow is induced by evaluation.  The result now follows
from \prettyref{lem:k}.
\end{proof}

\begin{cor}\label{cor:ineq2}
Let $G$ be a graph with at least one edge.  Then
\[\coind\Hom(K_2, G)\le\lim_{r\to\infty}\coind\Hom(C_{2r+1}, G)+1.\]
\end{cor}
\begin{proof}
Assume $\coind\Hom(K_2, G)\ge k+1$, i.e. the existence of a $\Z_2$-map
$\Sphere^{k+1}\to_{\Z_2}\Hom(K_2, G)$.  There is a map
$\Sphere^1_b\times_{\Z_2}\Sphere^k\to_{\Z_2}\Sphere^{k+1}$ and hence a
map $\Sphere^1_b\times_{\Z_2}\Sphere^k\to_{\Z_2}\Hom(K_2, G)$.  This
means that there is a map
\[\Sphere^k\to_{\Z_2}\Map_{\Z_2}(\Sphere^1_b,\Hom(K_2, G))
\homot_{\Z_2}\colim_{r}\Hom(C_{2r+1}, G).\] Because of the compactness of
$\Sphere^k$, it follows that whenever $r$ is large enough there is
a map $\Sphere^k\to_{\Z_2}\Hom(C_{2r+1}, G)$,
i.e.\ $\coind\Hom(C_{2r+1}, G)\ge k$.
\end{proof}

\subsection{Outline}
We start with necessary definitions and facts in
\prettyref{sec:prelims}.

In \prettyref{sec:proof} we give a short proof of
\prettyref{cor:ineq1}.  This reduces \prettyref{thm:S} to
\prettyref{thm:Lext}.  So far, the proof of \prettyref{thm:S} in
\cite{codd} had also been the simplest known proof of
\prettyref{thm:BK}.  However, \v Zivaljevi\'c had given a surprising
proof of
\[
\coind\Hom(C_{2r+1}, G)\le 2\left\lceil\frac{\chi(G)}2\right\rceil -3
\]
that was even simpler~\cite{parallel-transport,z-groupoids}.  Our
proof in \prettyref{sec:proof} extends arguments from \v
Zivaljevi\'c's proof. If one is content to obtain a bound on the
coindex of $\Hom(C_{2r+1}, G)$ instead of the cohomological index, it
is completely elementary in the sense that the Algebraic Topology used
is not more advanced than the degree of maps between spheres of the
same dimension, i.e.\ nothing more advanced than the Borsuk-Ulam Theorem.

In \prettyref{sec:gen} we generalize the proof and obtain conditions
under which inequalities similar to \prettyref{cor:ineq1} hold.  This
culminates in \prettyref{thm:k} which contains \prettyref{cor:ineq1}
as a special case.  As an application we show in \prettyref{ex:kneser}
that there is an infinite family of Kneser graphs with chromatic
number~$4$ which are test graphs.  Except for a topological lemma
at its beginning, this section is not needed for what follows.

In \prettyref{sec:colimit} we first prove \prettyref{cor:ineq2} and
finally \prettyref{thm:main-intro}.

\begin{acknowledgements}
I thank Mark de Longueville, Elmar Vogt and Rade \v Zivaljevi\'c for
helpful discussions.
\end{acknowledgements}

\section{Objects of study}\label{sec:prelims}
%%%%%%%%%%%%%%%%%%%%%%%%%%
We introduce the objects and concepts used in the proof.  The only
thing worth to be mentioned specifically is \prettyref{def:ast}, which
is very natural but to our knowledge has not been used explicitly
before.
\subsection{Free $\Z_2$-spaces}
We assume all spaces to be CW-spaces.
A good introduction to equivariant methods from the point of view of
combinatorial applications is~\cite{matousek}.
\begin{defn}\label{def:conn}
For an integer $m\ge-1$, we say that a topological space $X$ is
$m$-connected if every continuous map $\Sphere^k\to X$ with $-1\le
k\le m$ can be extended to a continuous map $\Disk^{k+1}\to X$.  We
define the \emph{connectivity} of $X$, $\conn X\in\Z\unite\set\infty$,
to be the largest~$m$ such that $X$ is $m$-connected.
\end{defn}

\begin{defn}\label{def:ind}
Let $X$ be a free $\Z_2$-space, i.e.\ a space with a fixed point free
involution.
We define the index and coindex of~$X$ by
\begin{align*}
\ind X&\deq\min\set{k\colon \text{There is a $\Z_2$-map $X\to\Sphere^k$}},
\\\coind X&\deq\max\set{k\colon \text{There is a $\Z_2$-map $\Sphere^k\to X$}}.
\intertext{Since there is a $\Z_2$-map $f\colon X\to\Sphere^\infty$ 
and this map is unique up to $\Z_2$-homotopy, we can also define
the cohomological index}
\hind X&\deq\max\set{k\colon \bar f^\ast(\gamma^k)\ne0},
\end{align*}
where $\bar f\colon X/\Z_2\to\RP^\infty$ is the map induced by~$f$,
and $H^\ast(\RP^\infty;\Z_2)=\Z_2[\gamma]$.
\end{defn}

\begin{defn}\label{def:even-odd}
If $X$, $Y$ are $\Z_2$-spaces and $f\colon X\to Y$ a map, then we will
also call $f$ an \emph{odd map} if it is equivariant, and we will call $f$
an \emph{even map} if it maps each orbit to a single point.
\end{defn}

\begin{prop}\label{prop:transfer}
If $X$ is a free $\Z_2$-space and $p_X\colon X\to X/\Z_2$ the
canonical quotient map, then there is the \emph{cohomology transfer
  maps} $p_X^!\colon H^\ast(X;\Z_2)\to H^\ast(X/\Z_2;\Z_2)$ which
fits in a long exact sequence
\begin{multline*}
H^k(X/\Z_2;\Z_2)\xto{p_X^\ast}H^k(X;\Z_2)\xto{p_X^!}H^k(X/\Z_2;\Z_2)
\\
\xto{\delta^\ast}H^{k+1}(X/\Z_2;\Z_2)
\xto{p_X^\ast}H^{k+1}(X;\Z_2)
\end{multline*}
which is natural with respect to $\Z_2$-maps.  For $X=\Sphere^\infty$
and $k\ge0$ it follows that $\delta^\ast\colon H^k(\RP^\infty;\Z_2)\to
H^{k+1}(\RP^\infty;\Z_2)$ is an isomorphism,
i.e.\ $\delta^\ast(\gamma^k)=\gamma^{k+1}$.
\end{prop}

Introductions to transfer maps which are sufficient for our purposes
when translated from homology to cohomology can be found in the proof
of the Borsuk-Ulam Theorem presented in the textbooks of
Bredon~\cite[pp.~240--241]{bredon} and Hatcher~\cite[p.~174]{hatcher}.
A more complete reference is \cite[Chap.~III]{bredon-groups}.

\subsection{Order complexes}
For a partially ordered set, or \emph{poset}, $P$, we denote by
$\Delta P$ its \emph{order complex}, the simplicial complex with
vertex set $P$ that consists of all chains in~$P$.  Any monotone (or
antitone) map $f\colon P\to Q$ between posets 
induces
a simplicial map
$\Delta f\colon\Delta P\to\Delta Q$. For a cell complex~$C$ we denote
its \emph{face poset} by $FC$ and its underlying space by $\real C$.
Thus $\Delta(FC)$ is the barycentric subdivision %\hfilneg\ \goodbreak 
of the complex~$C$
and $\real{\Delta(FC)}\homeo\real C$.  For posets $P$ and~$Q$,
$\Delta(P\times Q)$ is a simplicial subdivision of the cell complex
$\Delta P\times\Delta Q$, one often used for the product of simplicial
complexes with vertex orderings.

For posets $P$ and $Q$ we denote by $\Mon(P, Q)$ the poset of all
order preserving maps from $P$ to~$Q$.

\subsection{Graph complexes}
We will be brief in our description of $\Hom$-complexes.  A good
introduction is contained in~\cite{kozlov-survey}.

All graphs that we consider are finite, simple, and without loops.
The vertex set of a graph~$G$ is denoted by $V(G)$, the set of edges
by $E(G)$. 

\begin{nota}
Let $n\in\N$.  $K_n$ denotes the complete graph on
$n$~vertices with vertex set $\set{0,\dots,n-1}$.  $C_n$ is the cycle of
length~$n$ with vertex set $\set{0,\dots,n-1}$.
\end{nota}

\begin{defn}
A \emph{graph homomorphism} from $G$ to $H$ is a function
$f\colon V(G)\to V(H)$ that respects the edge relation, i.e.\ such
that $\set{f(u),f(u')}\in E(H)$ whenever $\set{u,u'}\in E(G)$.
The set of all graph homomorphisms from $G$ to~$H$ is denoted by
$\Hom_0(G, H)$.
\end{defn}

\begin{defn}
Let $G$, $H$ be graphs.  A \emph{multi-homomorphism} from $G$ to~$H$
is a function $\phi\colon V(G)\to\power(V(H))\wo\set\emptyset$ such
that every function $f\colon V(G)\to V(H)$ with $f(v)\in\phi(v)$ for
all~$v\in V(G)$ is a graph homomorphism.
\end{defn}

\begin{defn}\label{def:Hom}\sloppy
Let $G$, $H$ be graphs.  A function $\phi\colon
V(G)\to\power(V(H))\wo\set\emptyset$ can be identified with a cell of
the cell complex
\[\prod_{v\in V(G)}\Delta^{\#V(H)-1}.\]
The subcomplex of all cells indexed by multi-homomorphisms is denoted
by $\Hom(G, H)$.  We identify elements of $F\Hom(G, H)$ with the
corresponding multi-homomorphisms and $\Hom_0(G, H)$ with the
$0$-skeleton of $\Hom(G, H)$.
\end{defn}

\begin{defprop}\label{def:ast}\sloppy
The monotone map
\begin{align*}
\ast\colon F\Hom(G, G')\times F\Hom(G', G'')&\lto F\Hom(G, G'')\\
(\phi\ast\rho)(v)&\ \deq\ \rho[\phi(v)]
\intertext{induces a continuous map}
\ast\colon \real{\Hom(G, G')}\times\real{\Hom(G', G'')}&\lto
\real{\Hom(G, G'')}.
\end{align*}
This map is associative and its restriction to
$\Hom_0(G,G')\times\Hom_0(G',G'')\to\Hom_0(G,G'')$ coincides with
composition of graph homomorphisms.  Hence its restrictions
$\Hom_0(G,G')\times\Hom(G',G'')\to\Hom(G, G')$ and
$\Hom(G,G')\times\Hom_0(G',G'')\to\Hom(G, G')$ make $\Hom$ into a
functor, contravariant in the first and covariant in the second
argument.
\end{defprop}

\begin{defprop}
If $H$ is a graph and $\alpha\in\Hom_0(H, H)$ satisfies
$\alpha^2=\id_G$ and $\alpha$ flips an edge of~$H$ (i.e.\ the edge is
invariant but not fixed under~$\alpha$), then for every graph $G$
\begin{align*}
\real{\Hom(H, G)}&\to\real{\Hom(H, G)}\\
x&\mapsto \alpha\ast x
\end{align*}
is a fixed point free involution.  It is in this way that we make
$\Hom(K_2, G)$ and $\Hom(C_n, G)$ into free $\Z_2$-spaces.
\end{defprop}

\section{Proof of \prettyref{cor:ineq1}}\label{sec:proof}
%%%%%%%%%%%%%%%%%%%%%%%%%%%%%%%%%%%%%%%
Let $G$ be a graph, $r\ge1$.  We consider the free $\Z_2$-actions on
$\Hom(K_2, G)$ and $\Hom(C_{2r+1}, G)$ induced by
$\alpha\in\Hom_0(K_2, K_2)$ and $\beta\in\Hom_0(C_{2r+1},C_{2r+1})$
with $\alpha(v)=1-v$ and $\beta(v)=2r-v$.

Every $x\in\real{\Hom(K_2,C_{2r+1})}$ induces a map
\begin{align*}
f_x\colon\real{\Hom(C_{2r+1}, G)}&\to\real{\Hom(K_2,G)},\\
y&\mapsto x\ast y.
\end{align*}
We define $x_0,x_1\in F\Hom(K_2,C_{2r+1})$ by $x_0(0)\deq\set r$,
$x_0(1)\deq\set{r-1,r+1}$, $x_1(0)=\set 0$, $x_1(1)=\set{2r}$.  Then
$x_0\ast\beta=x_0$ and $x_1\ast\beta=\alpha\ast x_1$.  Hence,
$f_{x_0}(\beta\ast y)=f_{x_0}(y)$ and $f_{x_1}(\beta\ast y)=\alpha\ast
f_{x_1}(y)$ for all~$y$, i.e., using the language of \prettyref{def:even-odd}, 
$f_{x_0}$ is even and $f_{x_1}$ odd.

It is easy to check that the complex $\Hom(K_2, C_{2r+1})$ is
homeomorphic to a $1$-sphere and in particular path-connected.  A path
$(x_t)$ from $x_0$ to $x_1$ gives a homotopy $f_{x_t}$, and hence
$f_{x_0}\homot f_{x_1}$.

The inequality $\hind\Hom(C_{2r+1},G)+1\le\hind\Hom(K_2,G)$
now is a consequence of the following Lemma.\qed

\begin{lem}\label{lem:1}
Let $X$, $Y$ be free $\Z_2$-spaces, $Y\ne\emptyset$, 
$f,g\colon X\to Y$ maps.  If $f$
is odd, $g$ even, and $f\homot g$, then $\hind X+1\le \hind Y$.
\end{lem}

Before proving the Lemma we quickly show how to obtain the weaker
inequality $\coind X+1\le\ind Y$.  If for some $k$ there is a
$\Z_2$-map $k\colon Y\to\Sphere^k$, then the existence of a $\Z_2$-map
$l\colon\Sphere^k\to X$ would give rise to an odd map $k\cmps f\cmps
l\colon\Sphere^k\to\Sphere^k$ and an even map $k\cmps g\cmps l$.
These maps would be homotopic, which contradicts that even maps
between spheres have even degree and odd maps have odd
degree~\cite[Prop.~2B.6]{hatcher}.

\begin{proof}
Let $k\ge 0$.  We show that $\hind Y\le k$ implies $\hind X\lt k$.
Assume that $\bar h^\ast(\gamma^{k+1})=0$
where $\bar h$ is induced by a $\Z_2$-map $h\colon
Y\to\Sphere^\infty$.  From \prettyref{prop:transfer} we obtain 
a commutative diagram
\begin{equation*}\xymatrix{
&&H^k(\RP^\infty;\Z_2)\ar[r]_-\isom^-{\delta^\ast}\ar[d]^{\bar h^\ast}
&H^{k+1}(\RP^\infty;\Z_2)\ar[d]^{\bar h^\ast}
\\
&H^k(Y;\Z_2)\ar[r]^-{p_Y^!}\ar[ld]_{(g')^\ast}\ar[d]^{f^\ast}
&H^k(Y/\Z_2;\Z_2)\ar[r]^-{\delta^\ast}\ar[d]^{\bar f^\ast}
&H^{k+1}(Y/\Z_2;\Z_2)
\\
H^k(X/\Z_2;\Z_2)\ar[r]^-{p_X^\ast}
&H^k(X;\Z_2)\ar[r]^-{p_X^!}
&H^k(X/\Z_2;\Z_2)
}\end{equation*}
with exact rows, where $g'$ is defined by $g=g'\cmps p_X$.  Now 
\[\delta^\ast(\bar h^\ast(\gamma^k))=\bar h^\ast(\delta^\ast(\gamma^k))
 =\bar h^\ast(\gamma^{k+1})=0.
\]
Therefore there is an $\eta\in H^k(Y;\Z_2)$ with $\bar h^\ast(\gamma^k)=p_Y^!(\eta)$ and
\[
(\bar h\cmps\bar f)^\ast(\gamma^k)
=\bar f^\ast(\bar h^\ast(\gamma^k))
=\bar f^\ast(p_Y^!(\eta))=p_X^!(f^\ast(\eta))
=p_X^!(p_X^\ast((g')^\ast(\eta)))=0
\]
follows.
\end{proof}

\section{Generalizations}
%%%%%%%%%%%%%%%%%%%%%%%%%
\label{sec:gen}

In \prettyref{sec:proof} we have shown that \[\hind\Hom(C_{2r+1},
H)+1\le\hind\Hom(K_2,H), \] and the proof used properties of
$\Hom(K_2, C_{2r+1})$.  We will generalize this and determine
properties of $\Hom(G, G')$ which imply that
\begin{equation*}
\hind\Hom(G',H)+k\le\hind\Hom(G, H)
\end{equation*}
for a suitable $k\ge1$ and all graphs~$H$, \prettyref{thm:k}.  As an
application we obtain an infinite family of test graphs with
chromatic number~$4$ in \prettyref{ex:kneser}.

For the complex $\Hom(K_2, C_{2r+1})$ we implicitly used the
$\Z_2$-operation induced by the involution of $C_{2r+1}$ as well as
that induced by the involution of $K_2$.  It is easily seen that
there is a homeomorphism $\Hom(K_2, C_{2r+1})\homeo\Sphere^1$ under
which the operation induced by the involution of $C_{2r+1}$ is
equivalent to a reflection by a line through the origin of
$\R^2\supset\Sphere^1$ and the operation induced by the involution
of~$K_2$ is equivalent to the antipodal action on~$\Sphere^1$.  
This leads us to the following definition.

\begin{nota}\label{nota:Skb}
We write $\Z_2$ multiplicatively as $\Z_2=\set{1,\tau}$.  We will
always consider $\Sphere^k$ to be equipped with the left $\Z_2$ action
given by the antipodal map
\begin{align*}
\tau\cdot(x_0,\dots,x_k)&\deq(-x_0,\dots,-x_k).
\intertext{When additionally equipped with the right $\Z_2$-action by 
the reflection}
(x_0,\dots,x_k)\cdot\tau&\deq(-x_0,x_1,\dots,x_k)
\end{align*}
we will write the sphere as $\Sphere^k_b$.  Since these two actions
commute, we can also see them as a single $\Z_2\times\Z_2$-action.
\end{nota}

The case $k=1$ of the following lemma is equivalent to \prettyref{lem:1}.

\begin{lem}\label{lem:k}
Let $X\ne\emptyset$ be a free $\Z_2$-space and $k\ge0$. Then
\[\Sphere^k_b\times_{\Z_2}X\deq(\Sphere^k_b\times X)/_{(s,\tau x)\eqrel(s\tau,x)}\]
is also a free $\Z_2$-space, and
\[\hind X+k\le\hind (S^k\times_{\Z_2}X).\]
\end{lem}

\begin{proof}
Since the left and right action on~$\Sphere^k_b$ commute, the left
action induces an action on \hbox{$\Sphere^k_b\times_{\Z_2}X$}, which is free,
because the left actions on $\Sphere^k_b$ and~$X$ are free.  Furthermore
\hbox{$\Sphere^0_b\times_{\Z_2}X\homeo_{\Z_2}X$}, and it will be sufficient to show
\[\hind(\Sphere^{k-1}_b\times_{\Z_2}X)+1\le\hind(\Sphere^k_b\times_{\Z_2}X)\]
for $k\ge1$.

We consider the map
\begin{align*}
f\colon\Sphere^{k-1}\times[0,1]&\to\Sphere^k\\
(s,t)&\mapsto\left(\cos(t\pi/2)s,\sin(t\pi/2)\right).
\end{align*}
This map satisfies $f(s\cdot\tau,t)=f(s,t)\cdot\tau$, $f(\tau\cdot
s,0)=\tau\cdot f(s,0)$, $f(\tau\cdot s,1)=f(s,1)$ for all
$s\in\Sphere^{k-1}$ and~$t\in[0,1]$.  It therefore induces a homotopy
\begin{align*}
(\Sphere^{k-1}_b\times_{\Z_2}X)\times[0,1]&\to\Sphere^k_b\times_{\Z_2}X\\
([(s,x)], t)&\mapsto [(f(s,t),x)]
\end{align*}
from an odd to an even map and \prettyref{lem:1} can be applied.
\end{proof}

\begin{nota}
If $G$ and~$G'$ are graphs, then we write the $\Z_2$-action on
$\Hom(G, G')$ induced by an involution of $G$ as multiplication from
the left with elements of $\Z_2$, and the $\Z_2$-action induced by an
involution of $G'$ as right multiplication.  Again, these two actions
commute because of the associativity of~$\ast$.
\end{nota}

\begin{thm}\label{thm:Z2Z2}
Let $G$, $G'$ be graphs with edge-flipping involutions.  If there
exists a $\Z_2\times\Z_2$-map $f\colon \Sphere^k_b\to\Hom(G, G')$, then
\[\hind\Hom(G',H)+k\le\hind\Hom(G, H)\] for all graphs~$H$ with 
$\Hom(G', H)\ne\emptyset$.
\end{thm}

\begin{proof}
The diagram
\[\xymatrix{
\Sphere^k_b\times\Hom(G', H)
\ar[r]^-{f\times\id}\ar[d]&
\Hom(G, G')\times\Hom(G',H)
\ar[r]^-\ast\ar[d]&
\Hom(G, H)\\
\Sphere^k_b\times_{\Z_2}\Hom(G', H)
\ar[r]_>{\Z_2}&
\Hom(G, G')\times_{\Z_2}\Hom(G',H)
\ar[ru]_>{\Z_2}
}\]
commutes.  The existence of the $\Z_2$-map making the right triangle
commutative follows from the associativity of~$\ast$.  The two
$\Z_2$-maps at the bottom of the diagram show that \[\hind
(\Sphere^k\times_{\Z_2}\Hom(G', H))\le\hind\Hom(G, H),\] and
\prettyref{lem:k} concludes the proof.
\end{proof}

\begin{example}\label{ex:CoddKn}
As mentioned above, we indeed have
$\Hom(K_2,C_{2r+1})\homeo_{\Z_2\times\Z_2}\Sphere^1_b$, so that
\prettyref{cor:ineq1} is a special case of \prettyref{thm:Z2Z2}.
Similarly, if we equip $K_n$, $n>2$, with the $\Z_2$-action flipping
$\set{0,1}$ and keeping the other vertices fixed, then
$\Hom(K_2,K_n)\homeo_{\Z_2\times\Z_2}\Sphere^{n-2}_b$ for $n\ge 0$.
This can be seen as follows.  If for a non-empty subset $A$ of~$V(K_n)$
we let $b_A$ denote the barycentre of the corresponding face of the
$(n-1)$-simplex, then mapping $\phi\in F\Hom(K_2, K_n)$ to
$\frac12b_{\phi(0)}+\frac12b_{\compl\phi(1)}$ we obtain a homeomorphism
from $\Hom(K_2,K_n)$ to the boundary of the $\mbox{(n-1)}$-simplex.  This
sends the left action on $\Hom(K_2,K_n)$ to the antipodal map.  The
right action on $\Hom(K_2, K_n)$ is sent to the map on the simplex
induced by exchanging the vertices $0$ and~$1$.  This is a reflection
by the affine subspace spanned by $\set{b_{\set{0,1}},b_2,\dots,b_{n-1}}$.

We will take up these examples again in \ref{ex:Codd} and~\ref{ex:Kn}.
\end{example}

We would like to have a version of \prettyref{thm:Z2Z2} with
conditions that are easier to check.  A $\Z_2\times\Z_2$-map
$\Sphere^k_b\to\Hom(G, G')$ maps the fixed point sets in $\Sphere^k_b$ of
subgroups of $\Z_2\times\Z_2$ to the corresponding fixed point sets in
$\Hom(G, G')$.  Since the left action of $\Z_2$ on $\Sphere^k_b$ is
free, the fixed point set of $(\tau,1)$ in $\Sphere^k_b$ is empty.  The
fixed point set in $\Hom(G, G')$ of $(\tau,\tau)$ corresponds to the
equivariant multi-homomorphisms from $G$ to~$G'$.

\begin{defn}
For $i\in\set{0,1}$ let $G_i$ be a graph equipped with a $\Z_2$-action
given by $\alpha_i\in\Hom_0(G_i,G_i)$ with $\alpha_i^2=\id$.  We
define $\Hom_{\Z_2}(G_0,G_1)$ to be the subspace of
$\real{\Hom(G_0,G_1)}$ consisting of all $x$ with $\alpha_0\ast
x=x\ast\alpha_1$.
\end{defn}

The fixed-point set in $\Hom(G, G')$ of $(1,\tau)$ surely contains all
the multi-homo\-morph\-isms whose image is contained in the induced
subgraph of~$G$ on all vertices which are fixed by the involution.
However, since we are dealing with multi-homomorphism and not only
homomorphisms, the fixed point set can be larger.  This leads us to
the following definition.

\begin{defn}
Let $G$ be a graph equipped with a $\Z_2$-action given by a
homomorphism \hbox{$\alpha\colon G\to G$} with $\alpha^2=\id$.  We
define a graph $G^{\Z_2}$ by
\begin{align*}
V(G^{\Z_2})&\deq\set{\set{u,\alpha(u)}\colon u\in V(G)}\\
E(G^{\Z_2})&\deq\set{\set{\set{u,\alpha(u)},\set{v,\alpha(v)}}
    \colon\text{$\set{u,v}\in E(G)$ and $\set{u,\alpha(v)}\in E(G)$}}
\end{align*}
We also define $\iota_G\in F\Hom(G^{\Z_2},G)$ by $\iota_G(v)\deq v$.
\end{defn}

The graph $G^{\Z_2}$ is determined by the following universal
property.  
\begin{prop}
Let $G$ be graph and $\alpha$ an involution on $G$.  $\iota_G\colon
G^{\Z_2}\to G$ is a multi-homomorphism with
$\iota_G\ast\alpha=\iota_G$.  If $H$ is a graph and $h\colon H\to G$
is a multi-homomorphism with $h\ast\alpha=h$, then there is a unique
multi-homomorphism $g\colon H\to G^{\Z_2}$ with $h=g\ast\iota_G$.
\qed
\end{prop}
The analogous property with multi-homomorphisms replaced by
homomorphisms is fulfilled by the inclusion of the induced subgraph
of~$G$ on all vertices which are fixed by the involution.

We can now state a more practical, but slightly weaker, form of
\prettyref{thm:Z2Z2}.
\begin{thm}\label{thm:k}
Let $G,G'$ be graphs with $\Z_2$-actions, the action on $G$ flipping
an edge, and $k\ge1$.  If 
\begin{itemize}
\item $\coind\Hom(G,{G'}^{\Z_2})\ge k-1$,
\item $\Hom_{\Z_2}(G, G')\ne\emptyset$, and
\item $\Hom(G, G')$ is $(k-1)$-connected,
\end{itemize}
then
\begin{equation*}
\hind\Hom(G',H)+k\le\hind\Hom(G, H)
\end{equation*}
for all graphs $H$ with $\Hom(G', H)\ne\emptyset$.
\end{thm}

\begin{proof}[Proof of \prettyref{thm:k}]
Since $\coind\Hom(G,{G'}^{\Z_2})\ge k-1$, there is a map
\[h\colon\Sphere^{k-1}\to_{\Z_2}\Hom(G,{G'}^{\Z_2})\xto{\Hom(G, \iota_{G'})}\Hom(G,G').\]
This map satisfies $h(-s)=\tau h(s)$ and $h(s)\tau=h(s)$ for all
$s\in\Sphere^{k-1}$.  We also choose
$y\in\Hom_{\Z_2}(G,G')\subset\real{\Hom(G,G')}$, that is $\tau
y\tau=y$.  Since $\Hom(G,G')$ is $(k-1)$-connected, and the only fixed
point of the action $s\mapsto -s$ on $\Disk^k$ is the origin, these
choices can be extended to a map $g\colon\Disk^k\to\Hom(G, G')$ with
$g(-s)=\tau g(s)\tau$ for all~$s\in\Disk^k$, $g(0)=y$, and
$g|_{\Sphere^{k-1}}=h$.  We now define a map
\begin{align*}
f\colon\Sphere^k_b&\to\Hom(G, G')\\
(x_0,\dots,x_k)&\mapsto
\begin{cases}
g(x_1,\dots,x_k),&x_0\ge0,\\
g(x_1,\dots,x_k)\cdot\tau,&x_0\le0.
\end{cases}
\end{align*}
This map commutes with the left and right actions, so that
\prettyref{thm:Z2Z2} is applicable.
\end{proof}

\begin{cor}\label{cor:T}
Let $T$ be a graph with a $\Z_2$-action that flips an edge and
$k\ge1$.  If $\coind\Hom(K_2, T^{\Z_2})\ge k-1$ and $\Hom(K_2, T)$ is
$(k-1)$-connected, then
\begin{equation*}
\hind\Hom(T, H)+k\le\hind\Hom(K_2, H)
\end{equation*}
for all graphs $H$ with $\Hom(T,H)\ne\emptyset$.
It follows that
\begin{equation*}
\chi(G)\ge\hind\Hom(T, G)+k+2
\end{equation*}
for all graphs $G$ with $\Hom(T,G)\ne\emptyset$.  
In particular, if $k=\chi(T)-2$ then $T$ is a
test graph as defined in the introduction.
\end{cor}
\begin{proof}
For the first equation we set $G=K_2$, $G'=T$ in the Theorem.  The
edge of~$T$ that is flipped by the action ensures that
$\Hom_{\Z_2}(K_2, T)\ne\emptyset$.  Now for a graph $G$ with $\Hom(T,
G)\ne\emptyset$, we set $H=K_{\chi(G)}$ to obtain
\begin{multline*}
\hind\Hom(T, G)+k\le
\hind\Hom(T,K_{\chi(G)})+k
\le\\\le
\hind\Hom(K_2,K_{\chi(G)})=\chi(G)-2
\end{multline*}
and hence the second equation.
\end{proof}

\begin{example}[Odd circuits]\label{ex:Codd}
In \prettyref{sec:proof} we have dealt with $C_{2r+1}$, $r\ge1$.  The
result could also have been achieved by applying \prettyref{cor:T}
with $k=1=\chi(C_{2r+1})-2$.  The multi-homomorphism $x_0\in
F\Hom(K_2,C_{2r+1})$ used there shows that
$\Hom(K_2,C_{2r+1}^{\Z_2})\ne\emptyset$.  Indeed, $C_{2r+1}^{\Z_2}$ has the
single edge $\set{\set{r-1,r+1},\set r}$.
\end{example}

\begin{example}[Complete graphs]\label{ex:Kn}
If we equip $K_n$, $n>2$ with the $\Z_2$-action flipping $\set{0,1}$
and keeping the other vertices fixed, then we can apply
\prettyref{cor:T} with $k=n-2$, since $K_n^{\Z_2}\isom K_{n-1}$.  That
$K_n$ is a test graph has been shown in~\cite{babson-kozlov-i}.
In contrast to the case of $C_{2r+1}$, it will probably not come as a
surprise that other complete graphs do not yield better bounds on
chromatic numbers than~$K_2$.

Instead of applying \prettyref{cor:T} with $k=n-2$ to obtain
\[\hind\Hom(K_n,H)+n-2\le\hind\Hom(K_2, H)\] we can also apply
\prettyref{thm:k} with $k=1$ to abtain the stronger result
\[\hind\Hom(K_n, H)+1\le\hind\Hom(K_{n-1}, H).\]
The necessary fact that $\Hom(K_{n-1}, K_n)$ is path-connected follows
by observing that the group of permutations of $V(K_n)$ is generated
by transpositions of the form~$(i, n-1)$.  More generally,
$\Hom(K_m,K_n)$ is homotopy equivalent to a wedge of $(n-m)$-spheres
for $n\ge m$~\cite{babson-kozlov-i}.
\end{example}

\begin{example}[Kneser graphs]\label{ex:kneser}
Let $KG_{n,l}$ denote the Kneser graph of $n$-element subsets of
$\set{0,\dots,2n+l-1}$.  Edges are pairs of disjoint sets.  It is the
result of \cite{lovasz} that $\Hom(K_2, KG_{n,l})$ is
$(l-1)$-connected and $\chi(KG_{n,l})=l+2$.

Let $r,s\ge1$.  On the set $\set{0,\dots,4r+2s-1}$ we consider the
permutation \[\sigma\colon i\mapsto 4r+2s-1-i.\] This induces a
$\Z_2$-action on $KG_{2r,2s}$ which flips the edge
\[\set{\set{0,\dots,2r-1},\set{2r+2s,\dots,4r+2s-1}}.\]  There is a
homomorphism $KG_{r,s}\to KG_{2r,2s}^{\Z_2}$ given by $M\mapsto
\set{M\unite\sigma[M]}$.  Therefore
$\coind\Hom(K_2, KG_{2r,2s}^{\Z_2})\ge\coind\Hom(K_2,KG_{r,s})=s$, and
\prettyref{cor:T} can be applied with $k=s+1$.

For $s=1$, this yields that $KG_{2r,2}$ is a test graph with chromatic
number~$4$.  For $r>1$, it is triangle-free.  So far, the only known
test graph with chromatic number~$4$ was $K_4$.  There are also good
candidates in~\cite{z-groupoids}; these are built from triangles.
\end{example}

\section{The colimits of $\Hom(C_{2r+1},G)$ and $\Hom(C_{2r},G)$ 
for $r\to\infty$}
\label{sec:colimit}

\sloppy%!!!!!!!!!!!!!!!!

We fix a graph~$G$.

\begin{defn}
Let $m\ge3$. We define a monotone map
\begin{align*}
\eta_m\colon\FHom(C_m, G)&\to_{\Z_2}\Mon_{\Z_2}(\FHom(K_2,C_m),\FHom(K_2,G)),\\
\eta_m(\phi)(\rho)&\deq\rho\ast\phi.
\end{align*}
\end{defn}

This map, or rather its adjoint 
\[\FHom(K_2,C_m)\times_{\Z_2}\FHom(C_m, G)\to_{\Z_2}\FHom(K_2,G),\] 
has been vital in
\prettyref{sec:proof}.  We will now define a map that will induce a
homotopy inverse in the limit.

\begin{defn}
Let $m\ge3$. We define a monotone map
\begin{align*}
\theta_m\colon\Mon_{\Z_2}(\FHom(K_2,C_m),\FHom(K_2,G))
&\to_{\Z_2}\FHom(C_{3m}, G),\\
\theta_m(f)(3k)&\deq f((\set k,\set{k-1}))_2,\\
\theta_m(f)(3k+1)&\deq f((\set k,\set{k-1,k+1}))_1,\\
\theta_m(f)(3k+2)&\deq f((\set k,\set{k+1}))_2.
\end{align*}
All sums on the right hand side are to be understood modulo~$m$.
Elements of $\FHom(K_2, G)$ are written as pairs of subsets of~$V(G)$.
\end{defn}

\begin{lem}
The map $\theta_m$ is a well-defined $\Z_2$-map.
\end{lem}

\begin{proof}
Since $f((\set k, \set{k-1}))\subset f((\set k,\set{k-1,k+1}))$, every
element of $\theta_m(f)(3k+1)$ is a neighbour of every element of
$\theta_m(f)(3k)$, and similarly of every element of
$\theta_m(f)(3k+2)$.  Since
$\theta_m(f)(3k+3)=\theta_m(f)(3(k+1))=f((\set{k+1},\set k))_2=f((\set
k, \set{k+1}))_1$, every element of $\theta_m(f)(3k+2)$ is a neighbour
of every element of $\theta_m(f)(3k+3)$, and this calculation also
covers the case of the vertices $3m-1$ and~$0$.  Therefore
$\theta_m(f)$ is actually a multi-homomorphism from $C_{3m}$ to~$G$.  To show that $\theta_m$ is equivariant, we calculate
\begin{align*}
\theta_m(\tau\cdot f)(3k)&=(\tau\cdot f)((\set k,\set{k-1}))_2
    = f((\set{m-1-k},\set{m-k}))_2
   \\&=\theta_m(f)(3(m-k-1)+2)=\theta_m(f)(3m-1-3k)
   \\&=(\tau\cdot\theta_m(f))(3k),
\\
\theta_m(\tau\cdot f)(3k+1)&=(\tau\cdot f)((\set k,\set{k-1,k+1}))_1
    = f((\set{m-1-k},\set{m-k,m-k-2}))_1
   \\&=\theta_m(f)(3(m-k-1)+1)=\theta_m(f)(3m-1-(3k+1))
   \\&=(\tau\cdot\theta_m(f))(3k+1),
\\
\theta_m(\tau\cdot f)(3k+2)&=(\tau\cdot f)((\set k,\set{k+1}))_2
    = f((\set{m-1-k},\set{m-2-k}))_2
   \\&=\theta_m(f)(3(m-k-1))=\theta_m(f)(3m-1-(3k+2))
   \\&=(\tau\cdot\theta_m(f))(3k+2),
\end{align*}
which completes the proof.
\end{proof}

\begin{prop}\label{prop:theta-cont}
Let $G$ be a graph, $X$ a compact triangulable free $\Z_2$-space and
\[f\colon\Sphere^1_b\times_{\Z_2}X\to_{\Z_2}\real{\Hom(K_2,G)}\] an
equivariant map.  Then there exists an $r\ge1$ and an equivariant map
\[g\colon X\to_{\Z_2}\real{\Hom(C_{2r+1}, G)}.\]
\end{prop}

\begin{rem}
Using notation to be introduced later, this proposition can be seen to
provide a poor man's version of a continuous map
$\Map_{\Z_2}(\Sphere^1_b,\Hom(K_2, G))\to_{\Z_2}\colim_r\Hom(C_{2r+1},
G)$.
\end{rem}

\begin{proof}[Proof of \prettyref{prop:theta-cont}]
The complexes $\Hom(K_2, C_{2r+1})$ are triangulations
of~$\Sphere^1_b$.  If we take $r$ large enough and $K$ a
$\Z_2$-invariant triangulation of $X$ that is fine enough, then there
exists a monotone map
\[\mathrm \FHom(K_2, C_{2r+1})\times \mathrm F K
\isom\mathrm F(\Hom(K_2, C_{2r+1})\times K)\to_{\Z_2}\FHom(K_2, G)\] 
which induces an approximation of $f$.  This is adjoint to a montone map
\[\mathrm F K\to_{\Z_2}\Mon_{\Z_2}(\FHom(K_2, C_{2r+1}),\FHom(K_2, G)).\]
Its composition with $\theta_{2r+1}$ induces the desired map~$g$.
\end{proof}

\begin{cor}Let $G$ be a graph.  Then
\[\coind\Hom(K_2, G)\le1+\lim_{r\to\infty}\coind\Hom(C_{2r+1}, G).\]
\end{cor}

\begin{proof}
Let $\coind\Hom(K_2, G)\ge k+1$.  Since
$\Sphere^1_b\times_{\Z_2}\Sphere^k$ is a free $(k+1)$-dimensional
$\Z_2$-space, there exists a map
$\Sphere^1_b\times_{\Z_2}\Sphere^k\to\Hom(K_2, G)$.  By the preceding
proposition, there exist an~$r$ and a map
$\Sphere^k\to_{\Z_2}\Hom(C_{2r+1}, G)$.
\end{proof}

We now set up a framework that will allow us to pass to the limit of
the spaces $\Hom(C_m, G)$.
\begin{defn}
We define monotone maps
\begin{align*}
i_m\colon\FHom(C_m, G)&\to\FHom(C_{3m}, G),\\
i_m(\phi)(3k)&\deq\phi(k-1),\\
i_m(\phi)(3k+1)&\deq\phi(k),\\
i_m(\phi)(3k+2)&\deq\phi(k+1)
\intertext{and}
j_m\colon\Mon_{\Z_2}(\FHom(K_2,C_m),\FHom(K_2,G))
    &\to\Mon_{\Z_2}(\FHom(K_2,C_{3m}),\FHom(K_2,G)),\\
j_m(f)((A,\set{3k+1}\unite B))&\deq f((\set{k-1,k+1},k)),\\
j_m(f)((A\unite\set{3k+1},B))&\deq f((k,\set{k-1,k+1})),\\
j_m(f)((\set{3k},\set{3k-1}))&\deq f((\set{k-1},\set k)),\\
j_m(f)((\set{3k-1},\set{3k}))&\deq f((\set k,\set{k-1})).
\end{align*}
\end{defn}

\begin{prop}\label{prop:eta-theta}
$\theta_m\cmps\eta_m=i_m$ and $\eta_{3m}\cmps\theta_m\le j_m$.
\end{prop}

\begin{proof}
The calculations
\begin{align*}
\theta_m(\eta_m(\phi))(3k)
    &=\eta_m(\phi)((\set k, \set{k-1}))_2=(\phi(k),\phi(k-1))_2
     =\phi(k-1)=i_m(\phi)(3k),
\\
\theta_m(\eta_m(\phi))(3k+1)
    &=\eta_m(\phi)((\set k, \set{k-1,k+1}))_1
     \\&=(\phi(k),\phi(k-1)\unite\phi(k+1))_1
     =\phi(k)=i_m(\phi)(3k+1),
\\
\theta_m(\eta_m(\phi))(3k+2)
    &=\eta_m(\phi)((\set k, \set{k+1}))_2=(\phi(k),\phi(k+1))_2
     =\phi(k+1)=i_m(\phi)(3k+2)
\end{align*}
prove $\theta_m\cmps\eta_m=i_m$.  Furthermore
\begin{align*}
\eta_{3m}(\theta_m(f))((\set{3k},\set{3k-1}))
  &=(\theta_m(f)(3k),\theta_m(3(k-1)+2))
    \\&=(f((\set{k},\set{k-1}))_2,f((\set{k-1},\set{k}))_2)
    \\&=(f((\set{k-1},\set{k}))_1,f((\set{k-1},\set{k}))_2)
    \\&=f((\set{k-1},\set{k}))=j_m(f)(\set{3k},\set{3k-1}).
\end{align*}
Since $\theta_m(f)(3k-1)\subset\theta_m(f)(3k+1)$ and
$\theta_m(f)(3k+3)\subset\theta_m(f)(3k+1)$, we have
\begin{multline*}
\eta_{3m}(\theta_m(f))((A,\set{3k+1}\unite B))
    \le(\theta_m(f)(3k)\unite\theta_m(f)(3k+2),\theta_m(f)(3k+1))
    \\=(f((\set k,\set{k-1}))_2\unite f((\set k,\set{k+1}))_2,
          f((\set k,\set{k-1,k+1}))_1)
    \\\le f((\set{k-1,k+1},\set k))=j_m(f)((A,\set{3k+1}\unite B)).
\end{multline*}
This shows $\eta_{3m}\cmps\theta_m\le j_m$.
\end{proof}

\begin{defn}
For $m\ge3$ we define a graph homomorphism
\begin{align*}
\kappa_m\colon C_{m+2}&\to C_m\\
0&\mapsto m-1,\\
i&\mapsto i-1,\qquad 1\le i\le m,\\
m+1&\mapsto 0.
\end{align*}
Using arithmetic modulo~$m$ on the right hand side, this can simply be
written as $\kappa_m(i)=i-1$, which makes it clear that this
homomorphism commutes with the involutions on~$C_{m+2}$ and~$C_m$.  
For a graph $G$, we use the induced continuous maps
$\Hom(\kappa_m, G)$ to define the colimits (direct limits)
\[\colim_{\text{$m$ odd}}\Hom(C_m, G)\qquad\text{and}\qquad\colim_{\text{$m$
    even}}\Hom(C_m, G).\]  These carry induced $\Z_2$-actions.
\end{defn}
The choice of the graph homomorphisms $\kappa_m$ is not of great
importance, as the following lemma shows.
\begin{lem}\label{lem:3m}
Let $G$ be a graph.
The colimit of the diagram of all
\[\real{\Hom(C_{3^{n+1}}, G)}\xto{\abs{\Delta i_{3^{n+1}}}}\real{\Hom(C_{3^{n+2}}, G)}\]
is $\Z_2$-homotopy equivalent to $\colim_{\text{$m$ odd}}\Hom(C_m, G)$, and 
the colimit of the diagram of all
\[\real{\Hom(C_{4\cdot3^n}, G)}
\xto{\abs{\Delta i_{4\cdot3^n}}}\real{\Hom(C_{4\cdot3^{n+1}}, G)}\]
is $\Z_2$-homotopy equivalent to $\colim_{\text{$m$ even}}\Hom(C_m, G)$.
\end{lem}

\begin{proof}
The map $i_m$ is induced by a graph homomorphism $\iota_m\colon
C_{3m}\to C_m$.  We first consider the case of odd~$m$.  It is easy to
check that there is a path in $\Hom_{\Z_2}(C_{3(2r+1)}, C_{2r+1})$
connecting $\iota_{2r+1}$ and
$\kappa_{6r+1}\cdots\kappa_{2r+3}\kappa_{2r+1}$.  This induces a
$\Z_2$-homotopy between $\Delta i_{2r+1}$ and
$\Hom(\kappa_{6r+1},G)\cdots\Hom(\kappa_{2r+1}, G)$ and hence between
the homotopy colimits of the corresponding diagrams.  Since these
diagrams consist of simplicial inclusion maps, hence cofibrations, the
natural maps from their homotopy colimits to their colimits are also
homotopy equivalences.

For even~$r$ we proceed similarly, using a path from
$\iota_{6r}\iota_{2r}$ to
$\kappa_{18r-2}\cdots\kappa_{2r+2}\kappa_{2r}$.  These graph
homomorphisms are used because they agree on the vertices $0$,
$9r-1$, $9r$ and $18r-1$.
\end{proof}

We are now ready to prove the main theorem.
\begin{thm}\label{thm:homot-colim}
Let $G$ be a graph.  Then
\begin{align*}
\colim_{\text{$m$ odd}}\Hom(C_m, G)
    &\homot_{\Z_2}\Map_{\Z_2}(\Sphere^1_b,\Hom(K_2,G)),
\intertext{where $\Z_2$ acts on the right hand side by the right action on 
$\Sphere^1_b$, and}
\colim_{\text{$m$ even}}\Hom(C_m, G)
    &\homot_{\Z_2}\Map(\Sphere^1_b,\Hom(K_2,G)),
\end{align*}
where $\Z_2$ acts on the right hand side by the right action on
$\Sphere^1_b$ and the action on $\Hom(K_2, G)$.
\end{thm}

\begin{proof}
We use the diagrams from \prettyref{lem:3m}.
By \prettyref{prop:eta-theta} the diagram
\begin{equation*}
\xymatrix{
&\\ \real{\Hom(C_{3m},G)} \ar[r]^-{\real{\Delta\eta_{3m}}} \ar[u]
&\real{\Delta\Mon_{\Z_2}(\FHom(K_2, C_{3m}),\FHom(K_2, G))}
  \ar[u]
\\ \real{\Hom(C_{m},G)} \ar[r]^-{\real{\Delta\eta_m}} \ar[u]^{\real{\Delta i_m}}
&\real{\Delta\Mon_{\Z_2}(\FHom(K_2, C_{m}),\FHom(K_2, G))}
  \ar[lu]_-{\real{\Delta\theta_m}} \ar[u]^{\real{\Delta j_m}} 
\\\ar[u]&\ar[u] 
}
\end{equation*}
commutes up to homotopy and therefore induces a homotopy equivalence
between the homotopy colimits of the columns.  Since both columns
consist of simplicial inclusion maps, which are cofibrations, the
homotopy colimits are homotopy equivalent to the colimits.  

The $1$-dimensional cell-complex $\Hom(K_2, C_{3m})$ can be obtained
from $\Hom(K_2, C_m)$ by dividing each $1$-cell into three $1$-cells.
There is a corresponding homeomorphism $\real{\Delta(\FHom(K_2,
  C_{3m}))}\xto\homeo\real{\Delta(\FHom(K_2, C_m))}$ and the map
$\real{\Delta j_m}$ is induced by a map homotopic to it.
Thus there is a natural map from the homotopy colimit of the right column to
the space $\Map_{\Z_2}(\Hom(K_2, C_3), \Hom(K_2, G))$ in the case of
odd~$m$ respectively $\Map_{\Z_2}(\Hom(K_2, C_4), \Hom(K_2, G))$ in
the case of even~$m$.  Using the technique of the proof of
\prettyref{prop:theta-cont} we see that these maps are weak homotopy
equivalences and hence homotopy equivalences.  

All these constructions
can be carried out in such a way that the homotopy equivalences are
$\Z_2$-maps between free $\Z_2$-spaces and hence $\Z_2$-homotopy
equivalences.  Finally, $\Hom(K_2, C_3)\homeo\Sphere^1_b$ and the
quotient of $\Hom(K_2, C_4)$ by the free left $\Z_2$-action is
homeomorphic to $\Sphere^1_b$ as a right $\Z_2$-space.
\end{proof}

\begin{rem}
In particular, $\colim_r\Hom(C_{2r}, K_{n+2})$ is homotopy equivalent
to the free loop space of the $n$-sphere, a well-studied space.  In
\cite{kozlov-cycles} the cohomology groups of $\Hom(C_m, K_{n+2})$ are
determined.  Together with an analysis of the maps induced in
cohomology by $i_m$ or $\kappa_m$ this yields an elementary
calculation of the cohomology groups of free loop spaces of spheres.
\end{rem}

\begin{rem}
Since $\Hom(K_2, C_m)$ with the $\Z_2$-action induced by the action on
$\K_2$ is homeomorphic to $\Sphere^1$ with the antipodal action for
odd~$m$ and homeomorphic to $\Sphere^0\times\Sphere^1$ with an action
exchanging the components for even~$m$, we obtain
\begin{align*}
\colim_r\real{\Hom(C_{2r}, C_m)}\homot\Map(\Sphere^1,\real{\Hom(K_2, C_m)})
&\homot\coprod_{\Z}\Sphere^1,\\
\colim_r\real{\Hom(C_{2r+1}, C_m)}\homot\Map_{\Z_2}(\Sphere^1,\real{\Hom(K_2, C_m)})
&\homot\begin{cases}
\coprod_{\Z}\Sphere^1,&\text{$m$ odd},\\
\emptyset,&\text{$m$ even.}
\end{cases}
\end{align*}
The homotopy types of the spaces $\Hom(C_s,C_m)$ have been determined
in~\cite{sonja-dmitry-cycles}.
\end{rem}

\begin{rem}
We have seen that $\colim_r\Hom(C_{2r+1},
K_{n+2})\homot\Map_{\Z_2}(\Sphere^1, \Sphere^n)$.  There is a
canonical map
\[V_{2,n+1}\deq\set{(x,y)\in\Sphere^{n}\colon \langle x,y\rangle=0}
\xto\quad \Map_{\Z_2}(\Sphere^1, \Sphere^{n})\] which maps $(x,y)$
to a loop following the great circle through $x$ and~$y$ at constant
speed, starting at~$x$ in the direction of~$y$.  Among the spaces
$\Hom(C_{2r+1}, K_{n+2})$, the space $\Hom(C_5, K_{n+2})$ is special,
because it is a manifold~\cite{csorba-lutz}.  It has been conjectured
by Csorba~\cite{csorba-thesis} and proven in \cite{c5} that
there are homeomorphisms
\[\Hom(C_5, K_n)\homeo V_{2,n+1}.\]
It should not be difficult to check that these maps can be arranged in a diagram
\[\xymatrix{
\Hom(C_5, K_{n+2})\ar[r]\ar[d]^\homeo
&\colim_r\Hom(C_{2r+1},K_{n+2})
\ar[d]^\homot
\\
V_{2,n+1}\ar[r]
&
\Map_{\Z_2}(\Sphere^1,\Sphere^n)
}\]
which commutes up to homotopy.
\end{rem}

\begin{cor}
Let $G$ be a graph.  If $\Hom(K_2, G)$ is $(k+1)$-connected, then the
spaces $\colim_r\Hom(C_{2r+1}, G)$ and $\colim_r\Hom(C_{2r+1}, G)$ are
$k$-connected.
\end{cor}

\begin{proof}
Let $X$ be a connected free $\Z_2$-space with non-degenerate basepoint
$x_0$.  $\Map_{\Z_2}(\Sphere^1, X)$ is homeomorphic to
$\set{f\in\Map(I, X)\colon f(1)=\tau f(0)}$. Evaluating at~$0$ makes
this into the total space of a fibration over $X$ with fibre
$\set{f\in\Map(I, X)\colon f(0)=x_0, f(1)=\tau x_0}$.  The fibre is a
fibre of the path fibration over~$X$ and hence homotopy equivalent to
the loop space $\Omega X$.  If now $\pi_k(X)\isom\pi_{k+1}(X)\isom0$,
then from the part
\[0\isom\pi_{k+1}(X)\isom\pi_k(\Omega X)\to\pi_k(\Map_{\Z_2}(\Sphere^1, X))
  \to \pi_k(X)\isom0
\]
of the exact homotopy sequence of the fibration it follows that
$\pi_k(\Map_{\Z_2}(\Sphere^1, X))\isom0$.  The space
$\Map(\Sphere^1,X)$ is also the total space of a fibration with base~$X$
and fibre~$\Omega X$, so the same conclusion holds.
\end{proof}

%\section{Some topological lemmas}
%%%%%%%%%%%%%%%%%%%%%%%%%%%%%%%%%

\bibliographystyle{cas-ea}
\bibliography{math,topology,combi}

\end{document}